\documentclass[a4,12pt,reqno]{amsart}%
\usepackage[english]{babel}
\usepackage[latin1]{inputenc}
\usepackage{graphics}
\usepackage{ulem}
\usepackage{hhline}
\usepackage{dsfont}
\usepackage{mathrsfs}
\usepackage{fancyhdr}
\usepackage{amsmath,amssymb}
\usepackage{rotating}
\usepackage[latin1]{inputenc}
\usepackage[T1]{fontenc}
\usepackage{fancybox}
\usepackage{color}
\usepackage{colortbl}
\usepackage{setspace}
\usepackage{enumerate}
\usepackage[nice]{nicefrac}
\usepackage{amsthm}
\usepackage{multicol}
\usepackage{pifont}
\usepackage[french]{minitoc}
\usepackage{latexsym,amsfonts}
\usepackage{amsthm}
\usepackage{varioref}
\usepackage{textcomp}
\usepackage{lmodern}
\usepackage{mathpazo}
\usepackage{euscript}
\usepackage[pdftex]{hyperref}
\numberwithin{equation}{section}
\makeatletter\@addtoreset{equation}{section}

\DeclareMathSymbol{\subsetneqq}{\mathbin}{AMSb}{36}
\newtheorem {theorem}{Theorem}[section]
\newtheorem {Mtheorem}[theorem]{Main Theorem}
\newtheorem {definition}[theorem]{Definition}
\newtheorem {lemma}[theorem]{Lemma}
\newtheorem {proposition}[theorem]{Proposition}
\newtheorem {remark}[theorem]{Remark}

\newtheorem {corollary}[theorem]{Corollary}

\newcommand{\C}{\mathbb C}
\newcommand{\R}{\mathbb R}
\newcommand{\Z}{\mathbb Z}

\newcommand{\Ker}{K}
\newcommand{\Cte}{\h_o(\Gamma)}

\newcommand{\F}{\mathbf{F}}
\newcommand{\fin}{\hfill  
$\square$}

\newcommand{\p}{{\mathscr{P}}}

\newcommand{\Ws}{{\mathscr{W}}}

\newcommand{\h}{{\mathscr{H}}}

\newcommand{\norm}[1]{\left\Vert#1\right\Vert}

\newcommand{\set}[1]{\left\{#1\right\}}
\newcommand{\scal}[1]{\left<#1\right>}

\begin{document}

\title{Series and integral representations of the Taylor coefficients of the Weierstrass sigma-function}
\author{Allal Ghanmi}
 \address{E.D.P. and Spectral Geometry,
          Laboratory of Analysis and Applications-URAC/03,
          Department of Mathematics, P.O. Box 1014,  Faculty of Sciences,
          Mohammed V University, Rabat, Morocco}
          \email{ag@fsr.ac.ma}
 \author{Youssef Hantout}
 \address{UMR-CNRS 8524,  UFR Math., USTL,  Cit\'e Scientifique,
              59655 Villeneuve d'Ascq Cedex, France}
          \email{youssef.hantout@math.univ-lille1.fr} 
 \author{Ahmed Intissar}
 \address{E.D.P. and Spectral Geometry,
          Laboratory of Analysis and Applications-URAC/03,
          Department of Mathematics, P.O. Box 1014,  Faculty of Sciences,
          Mohammed V University, Rabat, Morocco}
          \email{intissar@fsr.ac.ma}
          
%

\maketitle
\begin{abstract}
 We provide two kinds of representations for the Taylor coefficients of the Weierstrass
$\sigma$-function $\sigma(\cdot;\Gamma)$ associated to an
arbitrary lattice $\Gamma$ in the complex plane $\C=\R^2$ - the
first one in terms of the so-called Hermite-Gauss series over
$\Gamma$ and the second one in terms of  Hermite-Gauss integrals
over $\C$.
\end{abstract}


\section{\bf Introduction and statement of main results}
\indent We let $\Gamma \subset \C$ be a given two-dimensional
lattice and denote by $\sigma(z;\Gamma)$ the associated
Weierstrass $\sigma$-function
 defined through the convergent infinite product
\begin{equation}\label{Weierstrass}
\sigma(z;\Gamma) = z\prod_{\gamma\in\Gamma\setminus\{0\}}
\Big(1-\frac z \gamma\Big)e^{\frac z \gamma + \frac 12 \left(\frac
z \gamma\right)^2}.
\end{equation}
 It is an odd entire function on $\C$ and therefore its Taylor series expansion at $z=0$ can be written as
\begin{equation}\label{TaylorW}
\sigma(z;\Gamma)=\sum_{r=0}^{\infty}
\Ws_r(\Gamma)\frac{z^{2r+1}}{(2r+1)!}.
\end{equation}
A representation of the Taylor coefficient $\Ws_r(\Gamma)$, in
terms of the Eisenstein series
\begin{equation} \label{EisensteinSeries}
 g_2=60  \sum_{\gamma\in\Gamma\setminus\{0\}} 
 \frac 1{\gamma^{4}} \quad \mbox{ and } \quad g_3 =140 \sum_{\gamma\in\Gamma\setminus\{0\}}    
 \frac 1{\gamma^{6}},
 \end{equation}
 has been obtained by K. Weierstrass \cite{Weierstrass} in 1882, to wit
\begin{equation}\label{TaylorWC}
\Ws_r(\Gamma) := \sum_{\begin{array}{c}  2m+3n=r  \\    m,n\geq 0
\end{array} }      a_{m,n}\left(\frac{g_2}2\right)^m(2g_3)^n,
\end{equation}
where the coefficients $a_{m,n}$ are such that $a_{0,0} =1$ and
$a_{m,n}=0$ whenever $m<0$ or $n<0$.
 Otherwise, they are connected by the recursion formula
\begin{align}
a_{m,n} =  3(m+1)a_{m+1,n-1} & + \frac{16}3(n+1)a_{m-2,n+1}  \label{RecForm}\\
 & -\frac 13(2m+3n-1)(4m+6n-1)a_{m-1,n} . \nonumber
\end{align}
The first few values of the coefficients $\Ws_r(\Gamma)$ in
\eqref{TaylorWC}; up to $r\leq 6$, are known to be given by
  \cite[page 391]{Magnus}
\begin{align*}\label{Coeff2}
 \sigma(z;\Gamma) & = z - \frac{g_2}{2^4\cdot 3\cdot 5 }z^5 - \frac{g_3}{2^3\cdot 3\cdot 5\cdot 7}z^7 -  \frac{g^2_2}{2^9\cdot 3^2\cdot 5\cdot 7}z^9
\\ & -
\frac{g_2g_3}{2^7\cdot 3^2\cdot 5^2\cdot 7\cdot 11} z^{11} +
\frac{(23g^3_2 - 576g^2_3)}{2^{10}\cdot 3^4\cdot 5^2\cdot 7\cdot
11\cdot 13}
  z^{13} + \cdots .\nonumber
\end{align*}

In the present paper, we shall give new representations of the
Taylor coefficients $\Ws_r(\Gamma)$ either as
series over $\Gamma$ of Hermite-Gauss type, or as integrals over $\C$ of the $\sigma$-function. 
Before doing so, we fix the notations we shall use repeatedly in
the sequel:
 Let $\set{\omega_1,\omega_2}$ be an oriented arbitrary $\R$-basis of $\C=\R^2$ and set $\Gamma=\Z\omega_1+\Z\omega_2$.
 The area of a fundamental cell $\Lambda(\Gamma)$ of $\Gamma$ is given by  $S=S_\Gamma=\Im(\overline{\omega_1}\omega_2)$, where $\Im(z)$
 denotes the imaginary part of $z\in\C$. We define $\nu \in \R^{>}:= ]0,+\infty[$ and $\mu \in \C$ by
\begin{equation}
\nu=\nu(\Gamma)=\frac\pi{S_\Gamma}
 \label{NuGamma}\end{equation}
and
\begin{equation} \label{MuGamma}
\mu =\mu(\Gamma)=\frac{i}{S_\Gamma}
\left({\zeta\left(\frac{\omega_1}{2};\Gamma\right)\overline\omega_2
-\zeta\left(\frac{\omega_2}{2};\Gamma\right)\overline\omega_1}\right)
.
\end{equation}
Here $\zeta(z;\Gamma)={\sigma^{'}(z;\Gamma)}/{\sigma(z;\Gamma)}$
is the Weierstrass zeta-function (not to be confused with the
Riemann zeta-function). Also, let $\chi_{_W}$ be the ``Weierstrass
semi-character" defined on $\Gamma$ by
\begin{equation}\label{Chi} \chi_{_W}(\gamma)= \left\{\begin{array}{ll}
                       +1 & \qquad \mbox{if} \quad  \gamma /2 \in \Gamma \\
                       -1 & \qquad \mbox{if} \quad  \gamma /2 \not\in \Gamma
                           \end{array}\right. .
\end{equation}
 Associated to the data of $\nu$, $\mu$ and $\chi_{_W}$, we consider the series
\begin{align}
 \h_r(\Gamma) = \nu  \frac{(2r+1)!}{2^r\cdot r!} \sum_{\gamma\in\Gamma}
              \mu^r{\F}\left(-r;\frac 32;-\frac{\nu^2\overline{\gamma}^2}{2\mu}\right)
               \chi_{_W}(\gamma)|\gamma|^2 e^{-\frac {\nu}2|\gamma|^2} \label{TaylorWM2}
 \end{align}
for every nonnegative integer $r$.
In \eqref{TaylorWM2}, ${\F}(a;c;x)$ denotes the usual confluent hypergeometric function 
and $|\gamma|^2$ the square of the Euclidean length of the vector
$\gamma\in\C=\R^2$. Note that \eqref{TaylorWM2} can be expressed
in terms of the Hermite polynomial $H_{2r+1}$ of degree $2r+1$ as
\begin{align}
 \h_r(\Gamma) =  \sum_{\gamma\in\Gamma}
                  \left(\sqrt{-\frac{\mu}{2}}\right)^{2r+1}H_{2r+1}\left(\frac{\nu \overline{\gamma}}{\sqrt{-2\mu}} \right)
               \chi_{_W}(\gamma) \gamma e^{-\frac {\nu}2|\gamma|^2}. \label{TaylorWM22}
  \end{align}

\begin{definition}
Since the quantity $\h_r(\Gamma)$, as given by \eqref{TaylorWM22},
involves the Hermite polynomial $H_{2r+1}$ as well as the Gaussian
function $e^{-\frac {\nu}2|\gamma|^2} $, we may call it
Hermite-Gauss series over $\Gamma$ of order $r$.
\end{definition}

With the above notations, the main results of the present paper
can be stated as follows:

\begin{Mtheorem}  \label{Thm1} Let $\Gamma$ be a lattice in $\C$, and  keep
$\nu=\nu(\Gamma)$, $\mu=\mu(\Gamma)$,
$\chi_{_W}=\chi_{_W}(\Gamma)$ and $\h_r(\Gamma)$ as above.  Then,
we have
\begin{itemize}
  \item[i)]  An expansion series over $\Gamma$ of the Weierstrass $\sigma$-function $\sigma(z;\Gamma)$:
  \begin{equation} \label{Identity10}
\sigma(z;\Gamma)= \frac{1}{\h_o(\Gamma)} e^{\frac \mu 2 z^2}
\sum_{\gamma\in\Gamma} \chi_{_W}(\gamma) \gamma e^{-\frac {\nu
}2|\gamma|^2+\nu  z\bar \gamma}.
\end{equation}

  \item[ii)] A reproducing integral formula over $\C$ for $\sigma(z;\Gamma)$:
\begin{align}\sigma(z;\Gamma)= \left(\frac{\nu}\pi\right) e^{\frac\mu 2 z^2} \int_{\C} e^{\nu z
\overline{w}-\frac\mu 2 w^2-\nu|w|^2} \sigma(w;\Gamma)
        dm(w), \label{RepKer25}
        \end{align}
        where $dm$ denotes the usual Lebesgue measure on $\C$.

 \item[iii)] A representation of the Taylor coefficient $\Ws_r(\Gamma)$ in terms of Hermite-Gauss series over $\Gamma$:
\begin{equation} \label{TaylorWM} \Ws_r(\Gamma) = \frac{\h_r(\Gamma)}{\h_o(\Gamma)} .\end{equation}

  \item[iv)] A representation of $\Ws_r(\Gamma)$ in terms of the Hermite-Gauss integral over $\C$:
\begin{align}\Ws_r(\Gamma) & = \left(\frac {\nu^2}\pi\right)  \frac{(2r+1)!}{2^r\cdot r!}  \label{RepKer26}\\
                        &\times \int_{\C} \overline{w}  \mu^r \, \F\left(-r;\frac 32;-\frac{\nu^2\overline{w}^2}{2\mu}\right)e^{-\frac\mu 2 w^2-\nu|w|^2} \sigma(w;\Gamma)dm(w) . \nonumber
\end{align}
\end{itemize}
\end{Mtheorem}

\begin{remark} It is likely impossible to proof directly that the quantity $$\h_o(\Gamma)=\nu\sum_{\gamma\in\Gamma} \chi_{_W}(\gamma) |\gamma|^2 e^{-\frac {\nu}2|\gamma|^2}$$ appearing in the denominator of \eqref{Identity10} and \eqref{TaylorWM}, is not zero. However,  we reassure the reader that $\h_o(\Gamma)\ne 0$ and a proof will be given in Section 4.
\end{remark}

\begin{remark}
The approach we will be adopting, to obtain representations
\eqref{Identity10} and \eqref{TaylorWM}, is based neither on the
partial differential equations satisfied by $\sigma(z;\Gamma)$ and
involving $g_2$ and $g_3$, as elaborated by Weierstrass in
\cite{Weierstrass},  nor on the bilinear operators as done in
\cite{Eilbeck}. In fact, we will derive them from a concrete
description of the
 functional space ${\mathcal{O}}^{\nu}_{\Gamma,\chi}(\C)$ of $(\Gamma,\chi)$-theta functions of magnitude $\nu$ on $\C$,
i.e., the space of all entire functions $f$ on $\C$ satisfying the
functional equation
\begin{equation} \label{FunEq}
f(z+\gamma):= \chi (\gamma) e^{\frac\nu 2|\gamma|^2 +\nu z
\overline{\gamma}}f(z).
\end{equation}
\end{remark}

As applications of Theorem \ref{Thm1}, one can obtain highly
nontrivial identities of arithmetic type
 and relate some modular forms to Hermite-Gauss series.  Moreover, full series of identities including those obtained by Perelomov \cite{Perelomov} can be deduced.
 Furthermore, our main results can be applied to show that, for every fixed integer $n\geq 2$, the Eisenstein series
 $G_{2n}(\Gamma)=\sum'1/\gamma^{2n}$ can be expressed as a rational fraction in the derivatives $T_j= \widetilde{\theta_{_W}}^{(2j-1)}(0;\Gamma)$; $j=1, \cdots,n$, where $\widetilde{\theta_{_W}}(z;\Gamma)$ is the Weierstrass-theta series
as defined by \eqref{Generator3}. The precise statements of the
above mentioned applications as well as their details will be
given in a separate paper.

The paper is outlined as follows. In Section 2, we review some
basic properties of the Hilbert space
${\mathcal{O}}^{\nu}_{\Gamma,\chi}(\C)$ of $(\Gamma;\chi)$-theta
functions. In Section 3, we introduce the modified Weierstrass
$\sigma$-function  and the Weierstrass $(\Gamma;\chi)$-theta
series, which have the important property of being generators of
${\mathcal{O}}^{\nu}_{\Gamma,\chi}(\C)$. Section 4 deals with the
proof of the main results summarized as Theorem \ref{Thm1}.


\section{\bf Background on $(\Gamma,\chi)$-theta entire functions}
In this section, we review some basic properties of the space
${\mathcal{O}}^{\nu}_{\Gamma,\chi}(\C)$ of $(\Gamma,\chi)$-theta
entire functions of magnitude $\nu>0$, i.e., the functional space
\begin{equation}\label{EspO}
{\mathcal{O}}^{\nu}_{\Gamma,\chi}(\C) =\set{f \,\,
\mbox{entire};\quad f(z+\gamma)= \chi(\gamma) e^{\frac{\nu}
2|\gamma|^2 + \nu z \overline{\gamma}}f(z), \, z\in \C, \gamma\in
\Gamma},
\end{equation}
where $\Gamma$ is a lattice in $\C$ and $\chi$ a given map on
$\Gamma$ with $|\chi(\gamma)|=1$ and satisfying the
({RDQ})-condition
$$\chi(\gamma+\gamma')= {\chi}(\gamma){\chi}(\gamma')e^{\frac{{\nu}} 2(\gamma\overline{\gamma'}
    -\overline{\gamma}\gamma')}, \qquad \gamma, \gamma' \in \Gamma. \eqno{\mbox{({RDQ})}}. $$

\begin{remark} \label{RemarkRDQC11}
By  $(RDQ)$ we refer to the ``Riemann-Dirac Quantization"
condition for the pair $(H,E)$, with $H(z,w)=(\nu/\pi)\scal{z,w}$
and $E(z,w)=\Im H(z,w)$.
\end{remark}

The $({RDQ})$ condition is in fact a necessary and sufficient
condition ensuring that ${\mathcal{O}}^{\nu}_{\Gamma,\chi}(\C)$ is
a nonzero vector space \cite{GI-JMP08}. Indeed, $(RDQ)$ is
equivalent to the condition that the complex valued
  function $J_{\nu,\chi}$ defined on $\Gamma\times\C^n$ by $J_{\nu,\chi}(\gamma,z):=\chi(\gamma )e^{i\nu E(z,\gamma)}$ is an
   automorphy factor satisfying the cocycle identity, $J_{\nu,\chi}(\gamma_1+\gamma_2,z)=J_{\nu,\chi}(\gamma_1,z+\gamma_2)J_{\nu,\chi}(\gamma_2,z)$.
Therefore, $\gamma(z;v) := (z+\gamma; \chi(\gamma )e^{\frac\nu
2|\gamma|^2+\nu z\overline{\gamma}} . v)$ defines an action of
$\Gamma$ on $\C^n\times\C$. The associated quotient space
$L:=(\C^n\times\C)/\Gamma$ is a holomorphic line bundle over the
torus $\C^n/\Gamma$ with fiber $\C=\tau^{-1}([z])$, where the
projection map $\tau: (\C^n\times\C)/\Gamma\longrightarrow
\C^n/\Gamma$ is the natural one induced from the canonical
projection $ \pi: \C^n\longrightarrow \C^n/\Gamma$. Whence, the
space ${\mathcal{O}}^{\nu}_{\Gamma,\chi}(\C)$ can be viewed as the
space of holomorphic sections of the holomorphic line bundle
$L=(\C\times\C)/\Gamma$. Its dimension is then given by the
Pfaffian  $\sqrt{\det E}$ of the associated skew-symmetric form
$E$ (see
\cite{Mumford74,Lang82,Frobenius1884,Bump02,PolishChuk02}). In our
case, $E(z,w):=(\nu/\pi)\Im \scal{z,w}$. Therefore, we have the
following

\begin{proposition} Let $\Gamma$ be a given lattice in $\C$, $\chi$ a map on $\Gamma$ with $|\chi(\gamma)|=1$ and $\nu>0$ such that the ({RDQ})-condition holds. Then, the functional space ${\mathcal{O}}^{\nu}_{\Gamma,\chi}(\C)$ is a finite dimensional space whose the dimension is given explicitly by
$$\dim {\mathcal{O}}^{{\nu}}_{\Gamma,{\chi}}(\C)=\left(\frac\nu \pi\right) S_\Gamma,$$
where $S_\Gamma$ is the area of a fundamental cell of $\Gamma$.
\end{proposition}

\begin{remark} Note that for given $f,g\in {\mathcal{O}}^{\nu}_{\Gamma,\chi}(\C)$, the function $z \mapsto f(z)\overline{g(z)}e^{-\nu|z|^2}$ is $\Gamma$-periodic. Thus, one can equip  ${\mathcal{O}}^{\nu}_{\Gamma,\chi}(\C)$ with the hermitian scalar product
\begin{equation}
\scal{\scal{f, g}}_\Gamma=\int_{\Lambda(\Gamma)} f(z)
\overline{g(z)} e^{-\nu|z|^2}dm(z)\label{InnerSP}
\end{equation}
which turns out to be a positive definite hermitian inner product
on ${\mathcal{O}}^{\nu}_{\Gamma,\chi}(\C)$. Thus the finite
dimensional Hilbert space ${\mathcal{O}}^{\nu}_{\Gamma,\chi}(\C)$
is a reproducing kernel Hilbert space. The next result gives an
explicit expression for its reproducing kernel.
\end{remark}

Consider the kernel function ${\Ker}^\nu_{\Gamma,\chi}(z,w)$  on
$\C\times \C$ given by the series
    \begin{equation}\label{RepKer2}
 {\Ker}^\nu_{\Gamma,\chi}(z,w) :=
 \left(\frac\nu \pi\right) e^{\nu z \overline{w} } \sum_{\gamma\in \Gamma}\chi(\gamma)e^{-\frac \nu 2|\gamma|^2+\nu\big( z \bar\gamma -\overline{w}\gamma \big)}.
\end{equation}
Then ${\Ker}^\nu_{\Gamma,\chi}(z,w)$ is well defined as a
holomorphic function in $z$ and antiholomorphic in $w$. Moreover,
the kernel function ${\Ker}^\nu_{\Gamma,\chi}(z,w)$ possesses the
properties summarized below.

\begin{theorem}
\label{LemPropO} Let ${\Ker}^\nu_{\Gamma,\chi}(z,w)$ be as in \eqref{RepKer2}. 
  \begin{itemize}
        \item[i)] For arbitrary fundamental cell $\Lambda(\Gamma) $  of the lattice $\Gamma$, we have
         \begin{equation}\label{Trace}\int_{\Lambda(\Gamma)} {\Ker}^\nu_{\Gamma,\chi}(z,z)
        e^{-\nu|z|^2}dm(z) = \left(\frac\nu \pi\right) S_\Gamma \ne 0.\end{equation}

 \item[ii)] The function ${\Ker}^\nu_{\Gamma,\chi}(z,w)$ satisfies the $\Gamma$-bi-invariance property
  \begin{equation}\label{GbiInv} {\Ker}^\nu_{\Gamma,\chi}(z+\gamma,w+\gamma') = \chi(\gamma) e^{\frac{\nu} 2|\gamma|^2 + \nu z\overline{\gamma}}
     {\Ker}^\nu_{\Gamma,\chi}(z,w) \overline{\chi(\gamma')} e^{\frac{\nu} 2|\gamma'|^2 + \nu\overline{w}\gamma'}\end{equation}
    for every $z,w\in \C$ and $\gamma, \gamma' \in \Gamma$.

 \item[iii)] Every $f\in {\mathcal{O}}^{\nu}_{\Gamma,\chi}(\C)$ can be reproduced either as
        \begin{equation}\label{Rep1} f(z) =\left(\frac\nu \pi\right) \int_{\C}  e^{\nu z \overline{w}}f(w)
        e^{-\nu|w|^2}dm(w) \end{equation}
        or as
         \begin{equation}\label{Rep2} f(z) =  \int_{\Lambda(\Gamma)}  {\Ker}^\nu_{\Gamma,\chi}(z,w) f(w)
        e^{-\nu|w|^2}dm(w) .\end{equation}
\end{itemize}
\end{theorem}

A sketched proof of the above theorem is presented below. However,
a detailed one can be found in \cite[Theorem 3.7]{GI-JMP08}.

\quad

\noindent {\it Sketched proof of Theorem \ref{LemPropO}.}
 i) follows by direct computation of
        $$\int_{\Lambda(\Gamma)} {\Ker}^\nu_{\Gamma,\chi}(z,z) e^{-\nu|z|^2}dm(z)$$
 using \eqref{RepKer2} and taking into account the $(RDQ)$ condition which implies
        $$\int_{\Lambda(\Gamma)} e^{\nu(z\bar \gamma-\bar z \gamma)}dm(z)=0$$
 for every nonzero $\gamma\in\Gamma.$

To prove ii), we begin by writing an expression for
${\Ker}^\nu_{\Gamma,\chi}(z+\gamma,w)$ using \eqref{RepKer2}.
Making a standard change of
     variable in the obtained series together with the $(RDQ)$ condition, we
     get
      $${\Ker}^\nu_{\Gamma,\chi}(z+\gamma,w) = \chi(\gamma)e^{\frac\nu 2|\gamma|^2+\nu z\bar \gamma}{\Ker}^\nu_{\Gamma,\chi}(z,w).$$
This gives rise to ii) since
${\Ker}^\nu_{\Gamma,\chi}(z,w)=\overline{{\Ker}^\nu_{\Gamma,\chi}(w,z)}$.

For a proof of iii), we note first that for every given $f\in
{\mathcal{O}}^{\nu}_{\Gamma,\chi}(\C)$, there exists a certain
constant $C>0$ such that $|f(z)|\leq Ce^{\frac\nu 2|z|^2}$. Then
by considering $$f_\varepsilon(z):=f(\varepsilon z)$$ for
$0<\varepsilon <1$, it follows that
  the holomorphic function $f_\varepsilon(z)$ satisfies the growth condition
$|f_\varepsilon(z)|\leq Ce^{\frac\nu 2\varepsilon^2|z|^2}$. Hence,
it belongs to the Bargmann-Fock space
$\mathcal{B}^{2,\nu}(\C)=\mathcal{O}(\C)\cap
L^2(\C;e^{-\nu|z|^2}dm(z)),$ whose the reproducing kernel is
$(\nu/\pi)e^{\nu z\overline{w}}$ \cite{Bargmann61}. Therefore, we
have
\begin{equation} \label{fVarepsilon}
f_\varepsilon(z) = \left(\frac\nu \pi\right)\int_{\C} e^{\nu
z\overline{w}}f_\varepsilon(w)e^{-\nu|w|^2}dm(w).
\end{equation}
Thus formula \eqref{Rep1} follows by letting $\varepsilon
\longrightarrow 1^{-}$ and applying the Lebesgue's dominated
convergence theorem. While \eqref{Rep2} can be obtained from
\eqref{Rep1} by writing $\C$ as a disjoint union of
$\gamma+\Lambda(\Gamma)$; $\gamma\in \Gamma$, and then using the
fact that $f$ belongs to ${\mathcal{O}}^{\nu}_{\Gamma,\chi}(\C)$.
\fin

\begin{corollary}\label{CorRepKer2} The function ${\Ker}^\nu_{\Gamma,\chi}(z,w)$ given by
 \eqref{RepKer2} is the reproducing kernel of the Hilbert space ${\mathcal{O}}^{\nu}_{\Gamma,\chi}(\C)$.
\end{corollary}

\begin{corollary}\label{CorRepKer00}
There exists at least one point $w_0\in \C$ such that the map
$z\mapsto {\Ker}^\nu_{\Gamma,\chi}(z,w_0)$ is not identically
zero.
\end{corollary}

\noindent {\it Proof.} In view of the bi-invariance property
\eqref{GbiInv}, it is clear that for every fixed $w\in\C$ the
entire function  $z\mapsto {\Ker}^\nu_{\Gamma,\chi}(z,w)$ belongs
to ${\mathcal{O}}^{\nu}_{\Gamma,\chi}(\C)$. The fact that
$z\mapsto {\Ker}^\nu_{\Gamma,\chi}(z,w_0)$ is not identically zero
for some $w_0\in \C$ follows from i) of Theorem \ref{LemPropO}.
\fin

\begin{remark}
According to the proof provided above, the existence of $w_0\in
\C$, for which $z\mapsto {\Ker}^\nu_{\Gamma,\chi}(z,w_0)$ is non
identically zero, is eventually not effective. Moreover, Remark
\ref{Rem6} below shows that $\omega_0$ belongs to $\C\setminus
\Gamma$.
\end{remark}

\section{\bf Generators of ${\mathcal{O}}^{\nu}_{\Gamma,\chi}(\C)$ associated to the Weierstrass semi-character}

From now on, we restrict ourselves to the particular case of the
triplet $(\Gamma, \nu, \chi)$, where $\Gamma$ is a fixed lattice
in $\C$, $\nu=\nu(\Gamma)= \pi/S_\Gamma$   and
$\chi=\chi_{_W}(\gamma)$ is the Weierstrass semi-character as in
\eqref{Chi},
\begin{equation}\label{Chi1} \chi=\chi_{_W}(\gamma):= \left\{\begin{array}{ll}
                       +1 & \qquad \mbox{if} \quad  \gamma /2 \in \Gamma \\
                       -1 & \qquad \mbox{if} \quad  \gamma /2 \not\in \Gamma
                           \end{array}\right. .
\end{equation}
 The following result is immediate (whose the proof is left to the reader):
\begin{lemma}\label{RDQ2W}  Let $\nu=\nu(\Gamma)= \pi/S_\Gamma$ as above. Then  $\chi=\chi_{_W}$ in \eqref{Chi1} satisfies the $(RDQ)$ condition
\begin{equation}\label{RDQ2}
\chi(\gamma+\gamma')= \chi(\gamma)\chi(\gamma')e^{\frac\nu
2(\gamma\overline{\gamma'}-\overline{\gamma}\gamma')}; \qquad
\gamma, \gamma' \in \Gamma.
\end{equation}
\end{lemma}

From Lemma \ref{RDQ2W} above and Theorem \ref{LemPropO}, we deduce
that ${\mathcal{O}}^{\nu}_{\Gamma,\chi}(\C)$ is a one dimensional
vector space. Two kinds of generators for
${\mathcal{O}}^{\nu}_{\Gamma,\chi}(\C)$ are constructed below -
the first one can be considered as a modified Weierstrass
$\sigma$-function, while the second one is a modified theta
function. Its construction is based on the data of the reproducing
kernel ${\Ker}^\nu_{\Gamma,\chi}(z,w)$ of
${\mathcal{O}}^{\nu}_{\Gamma,\chi}(\C)$.

\subsection{The modified Weierstrass $\sigma$-function.}
The main result of this section is the following

\begin{theorem}\label{KeyThm1}
There exist a unique real number $\nu>0$ and a unique complex
number $\mu$ such that the entire function
\begin{equation}\label{ModifiedSigma}
\widetilde{\sigma}_\mu(z;\Gamma):= e^{-\frac 12 \mu
z^2}\sigma(z;\Gamma); \quad z\in \C,
\end{equation}
is a generator of ${\mathcal{O}}^{\nu}_{\Gamma,\chi}(\C).$ More
precisely, $\nu$ and $\mu$ are explicitly given by
\begin{equation}\label{MuGammathm}
 \nu = \frac\pi {S_\Gamma} \qquad \mbox{and} \qquad \mu =  \frac i{S_\Gamma}\Big({ \zeta(\omega_1/2)\overline\omega_2
-\zeta(\omega_2/2)\overline\omega_1}\Big).
\end{equation}
\end{theorem}

\begin{definition}
The function
$\widetilde{\sigma}(z;\Gamma)=\widetilde{\sigma}_\mu(z;\Gamma)$
defined by (\ref{ModifiedSigma}) and associated to the special
value of $\mu$ given through \eqref{MuGammathm} will be called the
modified Weierstrass $\sigma$-function.
\end{definition}

\noindent {\it Proof of Theorem \ref{KeyThm1}.}
 We have to look for necessary and sufficient conditions on the pair $(\nu,\mu)$ so that the function $\widetilde{\sigma}_\mu(z;\Gamma)$, as given by
  \eqref{ModifiedSigma}, satisfies the functional equation
\begin{equation}\label{FuncEq2}
\widetilde{\sigma}_\mu(z+\gamma;\Gamma):= \chi(\gamma) e^{\frac\nu
2|\gamma|^2 +\nu z
\overline{\gamma}}\widetilde{\sigma}_\mu(z;\Gamma) ; \qquad
z\in\C, \, \gamma\in \Gamma .
\end{equation}
For this, we begin by recalling the pseudo-periodicity
\cite{Lang82,Jones-Singerman,PolishChuk02}, satisfied by the
standard Weierstrass $\sigma$-function,
\begin{equation*}
\sigma(z+\gamma;\Gamma) = \chi_{_W}(\gamma)
e^{(z+\gamma/2)\eta}\sigma(z;\Gamma),
\end{equation*}
where $\eta(\gamma)$ is defined by
$\eta(m\omega_1+n\omega_2)=m\eta_1+n\eta_2$, where
$\{\omega_1,\omega_2\}$ is a given oriented $\R$-basis of the
lattice $\Gamma=\Z\omega_1+\Z\omega_2$ and $\eta_j$; $j=1,2$, are
related to the Weierstrass zeta-function by
$\eta_j=2\zeta(\omega_j/2)$. Therefore,
$\widetilde{\sigma}_\mu(z;\Gamma)$ satisfies
 (\ref{FuncEq2}) if and only if $e^{ (\mu\gamma+\nu\overline{\gamma} -\eta)(z+\gamma/2)}=1$
for every $z\in\C$ and every $\gamma\in \Gamma$. It follows then
that
 the numbers $\nu$ and $\mu$ verify the linear system
\begin{equation*}\label{LinSystEq}
    \left(\begin{array}{cc}
                       \overline{\omega}_1 &  \omega_1  \\
                       \overline{\omega}_2 &  \omega_2
     \end{array}\right)
      \left(\begin{array}{c}
                       \nu \\
                        \mu
     \end{array}\right) =       \left(\begin{array}{c}
                       \eta_1 \\
                        \eta_2
     \end{array}\right)
\end{equation*}
whose the determinant $\overline{\omega}_1\omega_2
-\omega_1\overline{\omega}_2=2i S_\Gamma$ is not $0$. Whence,
\begin{equation*}\label{LinSystEq1}
           \left(\begin{array}{c}
                       \nu \\
                        \mu
     \end{array}\right) =  \frac 1{2iS}\left(\begin{array}{cc}
                       \omega_2 &  -\omega_1  \\
                       -\overline{\omega}_2 &  \overline{\omega}_1
     \end{array}\right)
      \left(\begin{array}{c}
                       \eta_1 \\
                        \eta_2
     \end{array}\right)    =  \frac 1{2iS}\left(\begin{array}{c}
                       \omega_2  \eta_1  - \omega_1 \eta_2  \\
                       -\overline{\omega}_2  \eta_1  +  \overline{\omega}_1 \eta_2
     \end{array}\right)
\end{equation*}
and therefore, the numbers $\nu$ and $\mu$ are explicitly given by
\begin{equation*}\label{MuGammaPf}
 \nu = \frac{\eta_1\omega_2    - \eta_2\omega_1 }{2iS_\Gamma} \qquad \mbox{and} \qquad
  \mu = \frac{\eta_2\overline\omega_1    - \eta_1\overline\omega_2 }{2iS_\Gamma} .
\end{equation*}
 The expression of $\nu$ reduces further to $ \nu = \pi /{S_\Gamma}$, making use of the Legendre's relation,
$\eta_1\omega_2 - \eta_2\omega_1 =2i\pi$ (\cite[page
102]{Jones-Singerman}). For such values of $\nu $ and $\mu$, the
nonzero function $\widetilde{\sigma}_\mu(z;\Gamma)$ belongs to the
one dimensional space ${\mathcal{O}}^{\nu}_{\Gamma,\chi}(\C)$ and
hence is a generator of ${\mathcal{O}}^{\nu}_{\Gamma,\chi}(\C)$.
\fin

\begin{remark}\label{Rem6}
Since ${\Ker}^\nu_{\Gamma,\chi}(z,w)$, in \eqref{RepKer2}, is the
reproducing kernel of the one dimensional Hilbert space $
\mathcal{O}^\nu_{\Gamma,\chi}(\C)$ (Corollary \ref{CorRepKer2}),
it follows that
$$
{\Ker}^\nu_{\Gamma,\chi}(z,w) =
\norm{\widetilde{\sigma}_\mu(\cdot,\Gamma)}^{-2}
\widetilde{\sigma}_\mu(z;\Gamma)
\overline{\widetilde{\sigma}_\mu(w;\Gamma)}
$$
with  $ \norm{\widetilde{\sigma}_\mu(\cdot,\Gamma)}^{-2}= - \nu^2
\h_o(\Gamma)/\pi$, so that $\h_o(\Gamma) <0$. Whence, the set of
zeros of $(z,w) \mapsto {\Ker}^\nu_{\Gamma,\chi}(z,w)$ is exactly
$(\Gamma \times \C)\cup (\C \times \Gamma)$.
\end{remark}

\subsection{The Weierstrass-theta series.}
Below, we consider another generator of
${\mathcal{O}}^{\nu}_{\Gamma,\chi}(\C)$ given as series over
$\Gamma$. Namely, we assert the following

\begin{theorem}\label{KeyThm2} Let $\nu=\pi/S_\Gamma$ and $\chi=\chi_{_W}$ as above. Then, the entire function
$\widetilde{\theta_{_W}}(z;\Gamma)$ defined by
\begin{equation}\label{Generator3}
\widetilde{\theta_{_W}}(z;\Gamma)=\sum\limits_{\gamma\in\Gamma}
\chi(\gamma) \gamma e^{-\frac {\nu }2|\gamma|^2+\nu  z\bar\gamma}
\end{equation}
is a generator of the space
${\mathcal{O}}^{\nu}_{\Gamma,\chi}(\C)$.
\end{theorem}

\begin{definition} The odd entire function
$\widetilde{\theta_{_W}}(z;\Gamma)$ given by \eqref{Generator3}
will be called the Weierstrass $(\Gamma,\chi)$-theta function.
\end{definition}

To prove Theorem \ref{KeyThm2}, we make use of the following
lemma.

 \begin{lemma}\label{PoincareLemma1} Let $\nu=\pi/S$ and $\chi=\chi_{_W}$ as above. Let $f$ be a holomorphic function on $\C$ such that
$|f(z)| \leq C e^{\alpha|z|^\beta}$ for some constant $C\geq 0$
and given real numbers $\alpha\geq 0$ and $\beta < 2$. Define
\begin{equation}\label{PoincareSeries}
[\p^{\nu}_{\Gamma,\chi}(f)](z):= \sum_{\gamma\in \Gamma}
\chi(\gamma)e^{-\frac\nu 2 |\gamma|^2 +\nu z\overline{\gamma}}
f(z-\gamma)
\end{equation}
to be the periodization \`a la Poincar\'e of the function $f$. Then,
 \begin{itemize}
    \item[i)] $\p^{\nu}_{\Gamma,\chi}(f)$ belongs to
    ${\mathcal{O}}^{\nu}_{\Gamma,\chi}(\C)$.
    \item[ii)] If in addition $f$ is an even function, then $\p^{\nu}_{\Gamma,\chi}(f)$ is also even and consequently
     is identically zero.
 \end{itemize}
\end{lemma}

\noindent {\it Proof of Theorem \ref{KeyThm2}.}
 By specifying $f(z)=1$ and appealing to ii) of Lemma \ref{PoincareLemma1}, we easily get the following identity,
 \begin{equation}\label{Identity00}
 \sum\limits_{\gamma\in\Gamma}\chi(\gamma) e^{-\frac {\nu }2|\gamma|^2+\nu z\overline{\gamma}} = 0
 \end{equation}
for every $z\in \C$. Hence, it follows that
\begin{equation}\label{Identity39}
\widetilde{\theta_{_W}}(z;\Gamma)
\stackrel{\eqref{Identity00}}{=}-\sum\limits_{\gamma\in\Gamma}\chi(\gamma)
(z-\gamma) e^{-\frac {\nu }2|\gamma|^2+\nu  z\bar \gamma}=
-[\p^{\nu}_{\Gamma,\chi}(z\mapsto z)](z),
\end{equation}
and therefore $\widetilde{\theta_{_W}}(z;\Gamma)$ belongs to
${\mathcal{O}}^{\nu}_{\Gamma,\chi}(\C)$, according to i) of Lemma
\ref{PoincareLemma1}.
 The hard part is to prove  that $\widetilde{\theta_{_W}}(z;\Gamma) $ is not identically zero on $\C$. To handle this, we recall from  Corollary
 \ref{CorRepKer00} that the holomorphic function
 $$  z \mapsto K^{\nu}_{\Gamma,\chi}(z,w_0) := \left(\frac\nu \pi\right)  e^{\nu z\overline{w}_0} \sum_{\gamma\in \Gamma}
    \chi(\gamma)e^{-\frac \nu 2|\gamma|^2+\nu (z\overline{\gamma}-\overline{w}_0\gamma)}, $$
    belonging to ${\mathcal{O}}^{\nu}_{\Gamma,\chi}(\C)$,
  is not identically zero on $\C$ for some $w_0$.
This means that $z \mapsto K^{\nu}_{\Gamma,\chi}(z,w_0)$ is a
generator of the one dimensional space
${\mathcal{O}}^{\nu}_{\Gamma,\chi}(\C)$.
 Since the modified Weierstrass function $\widetilde{\sigma}(z;\Gamma)$ is also a generator of ${\mathcal{O}}^{\nu}_{\Gamma,\chi}(\C)$, according to Theorem \ref{KeyThm1}, there exists a constant $C_{w_0}\ne 0$ such that
$K^{\nu}_{\Gamma,\chi}(z,w_0) =C_{w_0}
\widetilde{\sigma}(z;\Gamma)$. More explicitly, we have
\begin{equation}\label{PfRepKer2}
 \left(\frac\nu \pi\right)  e^{\nu z\overline{w}_0} \, \sum_{\gamma\in \Gamma}
    \chi(\gamma)e^{-\frac \nu 2|\gamma|^2+\nu (z\overline{\gamma}-\overline{w}_0\gamma)} = C_{w_0} e^{-\frac 12 \mu
z^2}\sigma(z;\Gamma).
\end{equation}
Differentiating both sides of \eqref{PfRepKer2} at $z=0$ and using
the fact $\sigma'(0;\Gamma)=1$ yield
\begin{equation*}
\left(\frac{\nu^2}\pi\right) \sum_{\gamma\in \Gamma}\chi(\gamma)
\overline{({w_0}+\gamma)} e^{-\frac \nu 2|\gamma|^2-\nu
\overline{w}_0\gamma} =C_{w_0}\ne 0.\end{equation*} Now, taking
the conjugate and changing $\gamma$ to $-\gamma$, keeping in mind
that the Weierstrass pseudo-character $\chi$ is real and even, we
obtain
\begin{equation*}
\left(\frac{\nu^2}\pi\right) \sum_{\gamma\in \Gamma}\chi(\gamma)
({w_0}-\gamma) e^{-\frac \nu 2|\gamma|^2+\nu w_0\overline{\gamma}}
=\overline{C_{w_0}}\ne 0.
\end{equation*}
Finally, by means of \eqref{Identity39}, we deduce that $
\widetilde{\theta_{_W}}(w_0;\Gamma)= - \frac\pi{\nu^2}
\overline{C_{w_0}} \ne 0$. This completes the proof of Theorem
\ref{KeyThm2}.  \fin

\quad

\noindent {\it Sketched proof of Lemma \ref{PoincareLemma1}.}
 The growth condition $|f(z)| \leq C e^{\alpha|z|^\beta}$, $\alpha\geq 0$, $\beta < 2$,
 satisfied by the holomorphic function $f$, ensures that the
Poincar\'e series \eqref{PoincareSeries}
 converges absolutely and uniformly on compact subsets of $\C$ and then is holomorphic.
To conclude i), note that for every given $\gamma_0\in \Gamma$, we
have
\begin{align*}
[\p^{\nu}_{\Gamma,\chi}(f)](z+\gamma_0)  := \sum_{\gamma\in
\Gamma} \chi(\gamma)e^{-\frac\nu 2 |\gamma|^2 +\nu (z+\gamma_0)
\overline{\gamma}} f(z+\gamma_0-\gamma).
 \end{align*}
 The change of summation index $\gamma=\gamma_0+\gamma'$ together with the $(RDQ)$ condition,
$$ \chi(\gamma_0+\gamma')= \chi(\gamma_0)\chi(\gamma')e^{\frac\nu 2(\gamma_0\overline{\gamma'}-\overline{\gamma_0}\gamma')},$$
yields
$$[\p^{\nu}_{\Gamma,\chi}(f)](z+\gamma_0)  = \chi(\gamma_0)e^{\frac \nu 2|\gamma_0|^2 + \nu z\overline{\gamma_0}}
 \sum_{\gamma'\in \Gamma} \chi(\gamma') e^{\nu  (\gamma_0\overline{\gamma'}-\overline{\gamma_0}\gamma')}
 e^{-\frac\nu2 |\gamma'|^2+\nu z\overline{\gamma'}} f(z-\gamma').$$
Now, since $e^{\nu
(\gamma_0\overline{\gamma'}-\overline{\gamma_0}\gamma')} =1$,
which is an immediate consequence of the $(RDQ)$ condition,  one
deduces
\begin{align*}
[\p^{\nu}_{\Gamma,\chi}(f)](z+\gamma_0)
& = \chi(\gamma_0)e^{\frac \nu 2|\gamma_0|^2 + \nu z\overline{\gamma_0}}\sum_{\gamma'\in \Gamma} \chi(\gamma') e^{-\frac\nu2 |\gamma'|^2+\nu z\overline{\gamma'}} f(z-\gamma')\\
&= \chi(\gamma_0)e^{\frac \nu 2|\gamma_0|^2 + \nu
z\overline{\gamma_0}}[\p^{\nu}_{\Gamma,\chi}(f)](z).
 \end{align*}
To get ii), we let $f$ be an even function satisfying the
hypothesis of Lemma \ref{PoincareLemma1}. Since $\chi$ is even, it
is easy to show
  that $\p^{\nu}_{\Gamma,\chi}(f)$ is also an even function belonging to the one dimensional space
${\mathcal{O}}^{\nu}_{\Gamma,\chi}(\C)$ which is generated by the
odd function $\widetilde{\sigma}(z;\Gamma)$. So necessarily
$\p^{\nu}_{\Gamma,\chi}(f)\equiv 0$ on $\C$.\fin

\quad

We conclude this section with the following remark.

\begin{remark} Assertion ii) of  Lemma \ref{PoincareLemma1} can be used to derive
  a full series of identities when specifying the function $f$. Here we provide some examples:

  \begin{enumerate}
  \item[i)] By considering the function $f(z)=\cos(\lambda z)$ with $\lambda\in\C$, one gets
$$ \sum_{\gamma\in\Gamma} \chi(\gamma)\cos(\lambda(z-\gamma))e^{-\frac{\nu }2|\gamma|^2+\nu  z\bar \gamma}=0,$$
 which for $z=0$ reduces to $\sum\limits_{\gamma\in\Gamma} \chi(\gamma)\cos(\lambda\gamma)e^{-\frac {\nu}2|\gamma|^2}=0.$

\item[ii)] Other interesting identities involve even polynomials.
Indeed, for $f(z)= z^{2p}$, $p=0,1,2, \cdots$, we obtain
\begin{eqnarray}
\sum_{\gamma\in\Gamma} \chi(\gamma)(z-\gamma)^{2p}e^{-\frac
{\nu}2|\gamma|^2+\nu z\bar \gamma}=0 .\label{Identity1}
\end{eqnarray}
In particular $\sum\limits_{\gamma\in\Gamma}\chi(\gamma)
\gamma^{2p} e^{-\frac {\nu}2|\gamma|^2}=0.$
 More generally, for every positive integer $k=0,1,2, \cdots,$ we have
\begin{eqnarray}
\psi_\Gamma(k):= \sum_{\gamma\in\Gamma}\chi(\gamma) \gamma^{k}
e^{-\frac {\nu }2|\gamma|^2}=0 . \label{Identity3}
\end{eqnarray}
Indeed, the odd case, i.e., $k=2p+1$ in \eqref{Identity3}, can be
handled using the change of the summation index taking into
account that $\chi$ is even.
 In fact this yields $\psi_\Gamma(k)=(-1)^k\psi_\Gamma(k)$, so that $\psi_\Gamma(2p+1)=0$.
\end{enumerate}
\end{remark}

\begin{remark}
The obtained identities (\ref{Identity3}), are exactly those
obtained by Perelomov \cite[Equation (47)]{Perelomov}, for regular
lattice with cell area $S =\pi$ (so $\nu=1$), when dealing with
the completeness of the coherent state system. The proof of the
simple case $p=0,$ given there requires a detailed knowledge of
the relationships between theta functions.
\end{remark}


In the next section, we proceed towards a proof of our main
result.

\section{\bf Proof of main results (Theorem \ref{Thm1}):}

\noindent {\bf Proof of i):} We have to prove the following
representation series  for the Weierstrass $\sigma$-function
\begin{equation} \label{RepW1}
\sigma(z;\Gamma)=  \frac{1}{\h_o(\Gamma)}e^{\frac\mu
2z^2}\sum_{\gamma\in\Gamma}\chi(\gamma) \gamma e^{-\frac {\nu
}2|\gamma|^2+\nu  z\bar \gamma},
\end{equation}
with $\h_o(\Gamma)=\nu
\sum_{\gamma\in\Gamma}\chi(\gamma)|\gamma|^2 e^{-\frac {\nu
}2|\gamma|^2}$. Indeed, since $\widetilde{\theta_{_W}}(z;\Gamma) $
and $\widetilde{\sigma}(z;\Gamma)$  are both generators of the one
dimensional space ${\mathcal{O}}^{\nu}_{\Gamma,\chi}(\C)$ (see
Theorem \ref{KeyThm2} and Theorem  \ref{KeyThm1}), there exists a
constant $C\ne 0$ such that $ \widetilde{\theta_{_W}}(z;\Gamma) =
C \widetilde{\sigma}(z;\Gamma)$. This reads explicitly as
\begin{equation} \label{Identity7}
\sum_{\gamma\in\Gamma}\chi(\gamma) \gamma e^{-\frac
{\nu}2|\gamma|^2+\nu  z\bar \gamma} =C  \,
\widetilde{\sigma}(z;\Gamma) = C \,  e^{-\frac 12 \mu
z^2}\sigma(z;\Gamma); \qquad C\ne 0 .
\end{equation}
Differentiation of both sides of (\ref{Identity7}) at $z=0$ yields
\begin{equation} \label{Identity8}
\nu \sum_{\gamma\in\Gamma}\chi(\gamma)|\gamma|^2 e^{-\frac
{\nu}2|\gamma|^2}= C \, \sigma'(0;\Gamma)= C \ne 0.
\end{equation}
This leads to $\h_o(\Gamma)=C\ne 0$, since the left hand side in
\eqref{Identity8} is exactly $\h_o(\Gamma)$. Using this fact in
\eqref{Identity7} gives rise to the representation series
\eqref{RepW1}.

\quad

\noindent {\bf Proof of ii):} According to Theorem \ref{KeyThm1},
the modified Weierstrass $\sigma$-function
$\widetilde{\sigma}(z;\Gamma):= e^{-\frac 12 \mu
z^2}\sigma(z;\Gamma)$ belongs to
${\mathcal{O}}^{\nu}_{\Gamma,\chi}(\C)$. By iii) of Theorem
\ref{LemPropO}, we get
\begin{align}\label{RepKer24}
\widetilde{\sigma}(z;\Gamma)= \left(\frac\nu \pi\right) \int_{\C}
e^{\nu z\overline{w}}
\widetilde{\sigma}(w;\Gamma)e^{-\nu|w|^2}dm(w)
        \end{align}
and therefore, we obtain the reproducing integral formula over
$\C$ for the Weierstrass $\sigma$-function
 \begin{equation} \label{RepKer241}
 \sigma(z;\Gamma)= \left(\frac\nu \pi\right)e^{\frac\mu 2 z^2} \int_{\C} e^{\nu z \overline{w}-\frac\mu 2 w^2}\sigma(w;\Gamma)e^{-\nu|w|^2}dm(w).
 \end{equation}
This completes the proof of ii).

\quad

 \noindent {\bf Proof of iii):} The expression for the Taylor coefficients of $\sigma(z;\Gamma)$ can be
 easily obtained by making use of  the obtained representation series
(\ref{RepW1})  combined with the following technical lemma.

\begin{lemma} \label{Identity101}
For arbitrary complex numbers $a$ and $b$, we have the following
expansion in power series of the exponential function $e^{a z^2+b
z}:$
\begin{align}\label{Expansion}
e^{az^2+bz} =  \sum_{r\geq 0} \frac{a^r}{r!}{\F}\left(-r;\frac
12;-\frac{b^2}{4a}\right) z^{2r}  +b \sum_{r\geq 0}
\frac{a^r}{r!}{\F}\left(-r;\frac 32;-\frac{b^2}{4a}\right)
z^{2r+1},
\end{align}
where ${\F}(a;c;x)$ is the usual confluent hypergeometric
function.
\end{lemma}
\noindent In fact, substitution of \eqref{Expansion} into the
representation series \eqref{RepW1}, with $a=\mu/2$ and
$b=\nu\overline{\gamma}$, gives
\begin{align*}\sigma(z;\Gamma) & =  \frac 1\Cte \sum_{r\geq 0} \frac{(\mu/2)^r}{r!}
\Big(\sum_{\gamma\in\Gamma} \chi(\gamma) \gamma {\F}\left(-r;\frac
12;-\frac{\nu^2\overline{\gamma}^2}{2\mu}\right) e^{-\frac {\nu
}2|\gamma|^2} \Big)z^{2r} \\ & +\frac \nu \Cte\sum_{r\geq 0}
\frac{(\mu/2)^r}{r!}\Big(\sum_{\gamma\in\Gamma}\chi(\gamma)|\gamma|^2{\F}\left(-r;\frac
32;-\frac{\nu^2\overline{\gamma}^2}{2\mu}\right)e^{-\frac {\nu
}2|\gamma|^2}  \Big) z^{2r+1}.\end{align*} The change of $\gamma$
by $-\gamma$ in the first summation of the right hand side and the
use of the fact
 $\chi(-\gamma)=\chi(\gamma)$ give rise to
\begin{equation*}
\sum_{\gamma\in\Gamma} \chi(\gamma) \gamma {\F}\left(-r;\frac
12;-\frac{\nu^2\overline{\gamma}^2}{2\mu}\right)e^{-\frac
{\nu}2|\gamma|^2}=0.
\end{equation*}
Therefore the above expression for $\sigma(z;\Gamma)$ reduces
further to the following one
\begin{equation*}
\sigma(z;\Gamma)= \sum_{r\geq 0} \frac{\h_r(\Gamma)}{\Cte}
\frac{z^{2r+1}}{(2r+1)!},\end{equation*} where we have set
\begin{equation*}
\h_r(\Gamma)= \nu\frac{(2r+1)!}{r!}
\left(\frac{\mu}2\right)^r\sum_{\gamma\in\Gamma}
\chi(\gamma)|\gamma|^2 {\F}\left(-r;\frac
32;-\frac{\nu^2\overline{\gamma}^2}{2\mu}\right)e^{-\frac {\nu
}2|\gamma|^2} .
\end{equation*}
This completes the proof of iii) of Theorem \ref{Thm1}.

\quad

\noindent {\bf Proof of iv):} This is an immediate consequence of
the reproducing integral formula \eqref{RepKer241} (i.e., ii) of
Theorem \ref{Thm1}).
 Indeed, by replacing the involved exponential function $e^{\frac\mu 2 z^2+\nu z\overline{w}}$
 by its expansion in power series as given by \eqref{Expansion} (with $a=\mu/2$ and $b=\nu
\overline{w}$), and comparing the obtained expression with
$$\sigma(z;\Gamma)=\sum_{r=0}^{\infty} \Ws_r(\Gamma)
\frac{z^{2r+1}}{(2r+1)!},$$ one easily deduces the representation
of the coefficients $\Ws_r(\Gamma)$ in terms of the Gauss-Hermite
integrals.

\quad

This completes the proof of Theorem \ref{Thm1}.  \fin

\begin{remark}
The technical Lemma \ref{Identity101} is basically an alternative
appropriate form of $e^{- z^2 + 2t z }$ which is the generating
function for the Hermite polynomials $H_{n}(t)$. Indeed, splitting
the obtained sum by collecting terms in $z^{2r}$ and those in
$z^{2r+1}$
 and then applying the transformations \cite[page 252]{Magnus}
             \begin{align*}
              H_{2r}(z) = (-1)^r\frac{(2r)!}{r!}  {\F}\left(-r;\frac 12;z^2\right)
             \end{align*}
 and
             \begin{align*}
            H_{2r+1}(z) =  (-1)^r\frac{(2r+1)!}{r!} \, 2z \, {\F}\left(-r;\frac 32;z^2\right)
             \end{align*}
 give rise to the result of Lemma \ref{Identity101}.

\end{remark}

\section*{Acknowledgements}
{ The authors are thankful to the anonymous
referee and to the editor  for their valuable suggestions for
improving the presentation of the paper.}

\end{document}